\documentclass[reqno,12pt]{amsart}

\usepackage{psfig}
\usepackage{amsmath,graphicx}
\usepackage{amssymb,amstext,amsthm}
\usepackage{latexsym}
\usepackage[norelsize]{algorithm2e}
\usepackage{amsfonts}
\usepackage{multirow}
\usepackage{array}
\usepackage{epstopdf}
\usepackage{bbm}
\usepackage{color}
\usepackage{float}
\usepackage{url}

\newtheorem{theorem}{Theorem}

\makeatletter
\def\BFGS{{\small BFGS}}
\def\LBFGS{{\small L-BFGS}}

\def\ones{\ensuremath{\mathbbm{1}}}

\newcommand{\minimize}[1]{{\displaystyle\minim_{#1}}}
\newcommand{\minim}{\mathop{\operator@font{minimize}}}
\newcommand{\mgap}{\;\;}

\newcommand{\defined}{\mathop{\,{\scriptstyle\stackrel{\triangle}{=}}}\,}
\newcommand{\subject}{\mathop{\operator@font{subject\ to}}}  
\newcommand{\words}[1]{\mgap\text{#1}\mgap}
\newcommand\diag{\mathop{\operator@font diag}\nolimits}

\renewcommand{\Re}{\mathbb{R}}

\begin{document}

\title{Trust-Region Methods for Sparse Relaxation}
%
\author[L. Adhikari]{Lasith Adhikari}
\email{ladhikari@ucmerced.edu}
\address{Applied Mathematics, University of California, Merced, Merced, CA 95343}

\author[J. Erway]{Jennifer B. Erway}
\email{erwayjb@wfu.edu}
\address{Department of Mathematics, Wake Forest University, Winston-Salem, NC 27109}

\author[S. Lockhart]{Shelby Lockhart}
\email{locksl12@wfu.edu}
\address{Department of Mathematics, Wake Forest University, Winston-Salem, NC 27109}

\author[R. Marcia] {Roummel F. Marcia}
\email{rmarcia@ucmerced.edu}
\address{Applied Mathematics, University of California, Merced, Merced, CA 95343}

\thanks{Lasith Adhikari's research is  supported by the UC Merced Graduate Student Opportunity Fellowship Program.}
\thanks{J.~B. Erway was supported by National
Science Foundation Grant  CMMI-1334042.}
\thanks{R.~F. Marcia is supported in part by National Science Foundation grant
CMMI-1333326.}

\date{\today}

%
%
%


\begin{abstract}
In this paper, we solve the $\ell_2$-$\ell_1$ sparse recovery problem by
transforming the objective function of this problem into an 
unconstrained differentiable function and apply a limited-memory
trust-region method.  Unlike gradient projection-type methods, which
uses only the current gradient, our approach uses gradients from previous
iterations to obtain a more accurate Hessian approximation.
Numerical experiments show that our proposed approach eliminates
spurious solutions more effectively while improving the computational
time to converge.
\end{abstract}

\maketitle

\section{Introduction}
\label{sec:intro}

This paper concerns solving the sparse recovery problem
\begin{equation}\label{eqn:LASSO}
	\minimize{f\in\Re^n}\mgap \frac{1}{2} \| Af - y \|_2^2 + \tau \| f \|_1,
\end{equation}
where $A \in \Re^{\tilde{m}  \times \tilde{n}}$, $f \in \Re^{\tilde{n}}$, $y \in \Re^{\tilde{m}}$, $\tilde{m} \ll \tilde{n}$, 
and $\tau > 0$ is a constant regularization parameter (see \cite{Tib96,CanT05,Don06}).  
By letting $f = u-v$, where $u, v \ge 0$, 
we write \eqref{eqn:LASSO} as the constrained
but differentiable optimization problem
\begin{eqnarray}
	\minimize{u, v\in\Re^{\tilde{n}}} \ & \ & \frac{1}{2} \| A(u-v) - y \|_2^2 + \tau \ones_n^T(u+v) \nonumber \\
	\text{subject to} && u, v \ge 0, 		\label{eqn:LASSOuv}
\end{eqnarray}
where $\ones_{\tilde{n}}$ is the $\tilde{n}$-vector of ones (see, e.g.,  \cite{Figueiredo2007}).
We transform \eqref{eqn:LASSOuv} into an unconstrained optimization problem 
by the change of variables 
$u_i = \log(1+e^{\tilde{u}_i})$ and $v_i = \log(1+e^{\tilde{v}_i}),$
where $\tilde{u}_i$, $\tilde{v}_i \in \Re$ for $1 \le i \le \tilde{n}$ (see \cite{BanMDG05,OhW15}).
With these definitions, $u$ and $v$ are guaranteed to be non-negative. 
Thus, \eqref{eqn:LASSOuv} is equivalent to the following minimization problem:
\begin{eqnarray}
	&& \hspace{-.6cm}  \underset{\tilde{u}, \tilde{v} \in\Re^{\tilde{n}}}{\text{min}}
	\Phi(\tilde{u}, \tilde{v})  \defined   \nonumber
	\frac{1}{2} \sum_{i=1}^{\tilde{m}} 
	\left [
	\left \{ \sum_{j=1}^{\tilde{n}}  A_{i,j} \log \left ( \frac{1+e^{\tilde{u}_j}}{1+e^{\tilde{v}_j}} \right ) \!  \right \} \! - \! y_i \!
	\right ]^2 \nonumber
	\\
	&& 
	\hspace{1.85cm}
	+ \
	\tau \sum_{j=1}^{\tilde{n}} \log \bigg (
	(1+e^{\tilde{u}_j})(1+e^{\tilde{v}_{j}})
	\bigg ).  \label{eqn:LASSOunc}
\end{eqnarray}
We propose 
solving \eqref{eqn:LASSOunc} using a limited-memory
quasi-Newton trust-region optimization approach, which we describe in
the next section.

\bigskip

\noindent \textbf{Related work.}  Quasi-Newton methods have been
previously shown to be effective for sparsity recovery problems (see
e.g., \cite{YuVGS10,LeeSS12,ZhoZD12}).  (For example, Becker and
Fadili use a zero-memory rank-one quasi-Newton approach for proximal
splitting \cite{BecF12}.)  Trust-region methods have also been
implemented for sparse reconstruction (see e.g.,
\cite{WanCY11,HinW13}).  Our approach is novel in the transformation
of the sparse recovery problem to a differentiable unconstrained
minimization problem and in the use of eigenvalues for efficiently
solving the trust-region subproblem.

\bigskip

\noindent \textbf{Notation.}  Throughout this paper, we denote the
identity matrix by $I$, with its dimension dependent on the context.

\section{Trust-Region Methods}
In this section, we outline the use of a trust-region method 
to solve \eqref{eqn:LASSOunc}.  We begin by combining the
unknowns $\tilde{u}$ and $\tilde{v}$ into one vector of unknowns
$x = [\tilde{u}^T \ \ \tilde{v}^T]^T \in \Re^{n}$, where
$n=2\tilde{n}$.  
(With this substitution,  $\Phi$ can be considered as a function
of $x$.)
Trust-region methods to minimize
$\Phi(x)$
define a sequence of iterates $\{x_k\}$
that are updated as follows:
$x_{k+1}=x_k+p_k,$
where $p_k$ is defined as the \emph{search direction}.
Each iteration, a new search direction $p_k$ is computed from
  solving the following 
quadratic subproblem with a two-norm constraint:
\begin{eqnarray} \label{eqn:TrustProblem}
  p_k &=& \underset{p \in\Re^n}{\text{arg min}} \hspace{.65cm} q_k(p) \defined g_k^Tp + \frac{1}{2} p^TB_kp \\
  	&& 
	\subject \mgap \|p\|_2 \le \delta_k,\nonumber
\end{eqnarray}
where $g_k\defined\nabla \Phi(x_k)$, $B_k$ is an approximation to
$\nabla^2 \Phi(x_k)$, and $\delta_k$ is a given positive constant.
In large-scale optimization, solving (\ref{eqn:TrustProblem}) represents
the bulk of the computational effort in trust-region methods.

Methods that solve the trust-region subproblem to high accuracy are
often based on the optimality conditions for a global solution
to the trust-region subproblem (see, e.g., \cite{Gay81,MorS83,ConGT00a})
given in the following theorem:

\medskip

\begin{theorem}\label{thrm-optimality}
  Let $\delta$ be a positive constant.  A vector $p^*$ is a global
  solution of the trust-region subproblem (\ref{eqn:TrustProblem}) if and only
  if $\|p^*\|_2\leq \delta$ and there exists a unique $\sigma^*\ge 0$ such
  that $B+\sigma^* I$ is positive semidefinite and
\begin{equation}\label{eqn-optimality}
(B+\sigma^* I)p^*=-g \mgap \words{and} \mgap \sigma^*(\delta-\|p^*\|_2)=0.
\end{equation}
Moreover, if $B+\sigma^* I$ is positive definite, then the global
minimizer is unique.
\end{theorem}


\section{Limited-Memory Quasi-Newton Matrices}
\label{sec:methods}
In this section we show how to build an approximation $B_k$ of
$\nabla^2 \Phi(x)$ using limited-memory quasi-Newton matrices.

Given the continuously differentiable function $\Phi$ 
and a sequence of iterates $\{x_k\}$, 
traditional quasi-Newton matrices are genererated from a sequence
of update pairs
$\{(s_k,y_k)\}$ where
$$
	s_k\defined x_{k+1} - x_k 
	\quad 
	\text{and}
	\quad 
	y_k \defined
	\nabla \Phi(x_{k+1}) - \nabla \Phi(x_{k}).
$$ 
In particular, given an initial matrix $B_0$, the
Broyden-Fletcher-Goldfarb-Shanno (\BFGS) update (see e.g.,
\cite{LiuN89,NocW06,GrNaS09}) generates a sequence of matrices
using the following recursion:
\begin{equation}\label{eqn:bfgs}
	B_{k+1}
	\defined 
	B_k 
	- \frac{1}{s_k^TB_ks_k}B_ks_ks_k^TB_k + \frac{1}{y_k^Ts_k}y_ky_k^T,
\end{equation}
provided $y_k^Ts_k \ne 0$.  In practice, $B_0$ is often taken to
be a nonzero constant multiple of the identity matrix, 
i.e., $B_0 = \gamma I$, for some $\gamma > 0$.  Limited-memory \BFGS{} 
(\LBFGS{}) methods store and use only the $m$ most-recently computed
pairs $\{(s_k,y_k)\}$, where $m\ll n$.  Often $m$ may be very small (for
example, Byrd et al.~\cite{ByrNS94} suggest $m\in [3,7]$).

The \BFGS{}
update is the most widely-used rank-two update formula that (i) satisfies
the \emph{secant condition} $B_{k+1}s_k=y_k$, (ii) has hereditary
symmetry, and (iii) provided that $y_i^Ts_i>0$ for $i=0,\ldots k$, then
$\{B_k\}$ exhibits hereditary positive-definiteness.


\bigskip 

\noindent \textbf{Compact representation.}
The \LBFGS{} matrix $B_{k+1}$ in \eqref{eqn:bfgs} can be defined recursively as follows:
$$
	B_{k+1}
	=
	B_0
	+ 
	\sum_{i=0}^{k} \left \{  - \frac{1}{s_i^TB_is_i}B_is_is_i^TB_i + \frac{1}{y_i^Ts_i}y_iy_i^T \right \}.
$$
Then $B_{k+1}$ is at most a rank-$2(k+1)$ perturbation to $B_0$, and 
thus, $B_{k+1}$ can be written as
$$
	B_{k+1} \ = \ B_0 + 
	\begin{bmatrix}
	\\
	\Psi_k  \\
	\phantom{t}
	\end{bmatrix}
	\hspace{-.3cm}
	\begin{array}{c}
	\left  [ \ M_k^{\phantom{h}}  \right ] \\
	\\
	\\
	\end{array}
	\hspace{-.3cm}
	\begin{array}{c}
	\left [  \ \quad \Psi_k^T \quad \ \right ] \\
	\\
	\\
	\end{array}
$$
for some $\Psi_k \in \Re^{n \times 2(k+1)}$ and $M_k \in \Re^{2(k+1)  \times 2(k+1)}$.
Byrd et al.~\cite{ByrNS94}
showed that $\Psi_k$ and $M_k$ are given by
\begin{equation*}
	\Psi_k \ = \ 
       	\begin{bmatrix}
		B_0 S_k  \ &  Y_k
	\end{bmatrix}
	\ \text{and} \ \
        M_k 
        \ = \ 
        -
        \begin{bmatrix}
        S_k^TB_0S_k & \ \ L_k \\
        L_k^T & -D_k
        \end{bmatrix}^{-1},
        \end{equation*}
where
\begin{eqnarray*}
	S_k &\defined& [ \ s_0 \ \ s_1 \ \ s_2 \ \ \cdots \ \ s_{k} \ ] \ \in \ \Re^{n \times (k+1)}, \\
	Y_k &\defined& [ \ y_0 \ \ y_1 \ \ y_2 \ \ \cdots \ \ y_{k} \ ] \ \in \ \Re^{n \times (k+1)},
\end{eqnarray*}
and $L_k$ is the strictly lower triangular part and $D_k$ is the diagonal part of the matrix
$S_k^TY_k \in\Re^{(k+1) \times (k+1)}$:
$$
	S_k^TY_k =   L_k + D_k + U_k.
$$
(In this decomposition, $U_k$ is a strictly upper triangular matrix.)


\section{Solving the Trust-region Subproblem}

In this section, we show how to solve \eqref{eqn:TrustProblem} efficiently.  
First, we transform \eqref{eqn:TrustProblem} into an equivalent expression.
For simplicity, we drop the subscript $k$.  Let $\Psi = QR$ be the ``thin'' QR
factorization of $\Psi$, where $Q \in \Re^{n \times 2(k+1)}$ has orthonormal columns
and $R \in \Re^{2(k+1) \times 2(k+1)}$ is upper triangular.  Then
$$
	B_{k+1} = B_0 + \Psi M \Psi^T = \gamma I + QRMR^TQ^T.
$$
Now let 
$
	V\widehat{\Lambda}V^T = RMR^T
$ 
be the eigendecomposition of 
$RMR^T \in \Re^{2(k+1)\times2(k+1)}$, where 
$V \in \Re^{2(k+1)\times2(k+1)}$ is orthogonal and 
$\widehat{\Lambda}$ 
is diagonal with $\widehat{\Lambda} = $ diag($\hat{\lambda}_1, \dots, \hat{\lambda}_{2(k+1)}$).
We assume that the eigenvalues $\hat{\lambda}_i$
are ordered in increasing values, i.e., $\hat{\lambda}_1 \le \hat{\lambda}_2 \le \cdots 
\le \hat{\lambda}_{2(k+1)}$.
Since $Q$ has orthonormal columns and $V$ is orthogonal, then $P_\parallel \defined
QV\in\Re^{n\times 2(k+1)}$ also has orthonormal columns.  Let $P_\perp$ be a
matrix whose columns form an orthonormal basis for the orthogonal
complement of the column space of $P_\parallel$.  Then, $P \defined [ \
P_{\parallel} \ \ \ P_{\perp} ] \in \Re^{n \times n}$ is such that $P^TP = PP^T =
I$.  Thus, the spectral decomposition of $B$ is given by
%
%
\begin{equation}\label{eqn-Beig}
	\phantom{.} \hspace{-.3cm} B = P\Lambda P^T, \ 
	\text{where }
	\Lambda  \defined
	\begin{bmatrix}
		\Lambda_1 & 0 \\
		0 & \Lambda_2
	\end{bmatrix} = 
	\begin{bmatrix}
		\hat{\Lambda} + \gamma I & 0 \\
		0 & \gamma I
	\end{bmatrix},
\end{equation}
where 
$\Lambda_{\phantom{1}} = \diag(\lambda_1,\ldots,\lambda_n),$
$\Lambda_1 = \diag(\lambda_1,\ldots, \lambda_{2(k+1)}) \in \Re^{2(k+1)\times 2(k+1)},$,
and $\Lambda_2 = \gamma I_{n-2(k+1)}$.
Since the $\hat{\lambda}_i$'s are ordered, then the eigenvalues in $\Lambda$ are also ordered, i.e.,
$\lambda_1\le \lambda_2 \le \ldots \le \lambda_{2(k+1)}$.  The remaining eigenvalues, found
on the diagonal of $\Lambda_2$, are equal to $\gamma$.  
Finally, since $B$ is positive definite, then $0 < \lambda_i$ for all $i$.

Defining $v = P^Tp$, the trust-region subproblem \eqref{eqn:TrustProblem}, 
can be written as
\begin{eqnarray} \label{eqn:TrustProblem2}
  v^* &=& \underset{v \in\Re^n}{\text{arg min}} \hspace{.65cm} q_k(v) \defined \tilde{g}^Tv + \frac{1}{2} v^T\Lambda v \\
  	&& 
	\subject \mgap \|v\|_2 \le \delta,\nonumber
\end{eqnarray}
where $\tilde{g} = P^Tg$.  
From the optimality conditions in Theorem 1,
the solution, $v^*$, to \eqref{eqn:TrustProblem2} must satisfy the following equations:
\begin{eqnarray}
	(\Lambda + \sigma^* I)v^* &=& -\tilde{g}  \\
	\sigma^* ( \| v^* \|_2 - \delta) &=& 0 \label{eqn:comp} \\
	\sigma^* &\ge& 0 \\
	\| v^* \|_2 &\le& \delta,
\end{eqnarray}
for some scalar $\sigma^*$.  Note that the usual requirement that
$\sigma^* + \lambda_i \ge 0$ for all $i$ is not necessary here since $\lambda_i > 0$
for all $i$ (i.e., $B$ is positive definite).
Note further that \eqref{eqn:comp}  implies that if $\sigma^* > 0$, 
the solution must lie on the boundary, i.e., $\| v^* \|_2 = \delta$.  In this case, the optimal 
$\sigma^*$ can be obtained by finding solving the so-called \emph{secular equation}:
\begin{equation}\label{eqn:secular}
	\phi(\sigma) = \frac{1}{\| v(\sigma) \|_2} - \frac{1}{\delta} = 0,
\end{equation}
where $\| v(\sigma) \|_2 = \| -(\Lambda + \sigma I)^{-1}\tilde{g} \|_2$.  
Since $\lambda_i + \sigma > 0$ for any $\sigma \ge 0$, 
$v(\sigma)$ is well-defined.  In particular, if we let
$$
	\tilde{g} = 
	\begin{bmatrix}
		P_{||}^T \\
		P_{\perp}^T
	\end{bmatrix}
	g 
	=
	\begin{bmatrix}
		P_{||}^Tg \\
		P_{\perp}^Tg
	\end{bmatrix}
	=
	\begin{bmatrix}
		g_{||} \\
		g_{\perp}
	\end{bmatrix},
$$
then
\begin{equation}\label{eq:normv}
	\| v(\sigma) \|_2^2 = \left \{ \sum_{i=1}^{2(k+1)} \frac{(g_{||})_i^2}{(\lambda_i - \sigma)^2} \right \} 
				+ \frac{\|g_{\perp}\|_2^2}{(\gamma - \sigma)^2}.
\end{equation}
We note that $\phi(\sigma) \ge 0$ means $v(\sigma)$ is feasible, i.e.,
$\| v(\sigma) \|_2 \le \delta$.  Specifically, the unconstrained
minimizer $v(0) = -\Lambda^{-1}\tilde{g}$ is feasible if and only if
$\phi(0) \ge 0$ (see Fig.\ 1(a)).  If $v(0)$ is not feasible, then
$\phi(0) < 0$ and there exists $\sigma^* > 0$ such that $v(\sigma^*) =
-(\Lambda + \sigma^* I)^{-1}\tilde{g}$ with $\phi(\sigma^*) = 0$ (see
Fig.\ 1(b)).  Since $B$ is positive definite, the function
  $\phi(\sigma)$ is strictly increasing and concave down for $\sigma \ge 0$, making it a
  good candidate for Newton's method.  In fact, it can be shown that
  Newton's method will converge monontonically and quadratically to
  $\sigma^*$ with initial guess $\sigma^{(0)}=0$~\cite{ConGT00a}.

\begin{figure}[h]
	\centering
	\begin{tabular}{cc}
	\includegraphics[width=4.2cm]{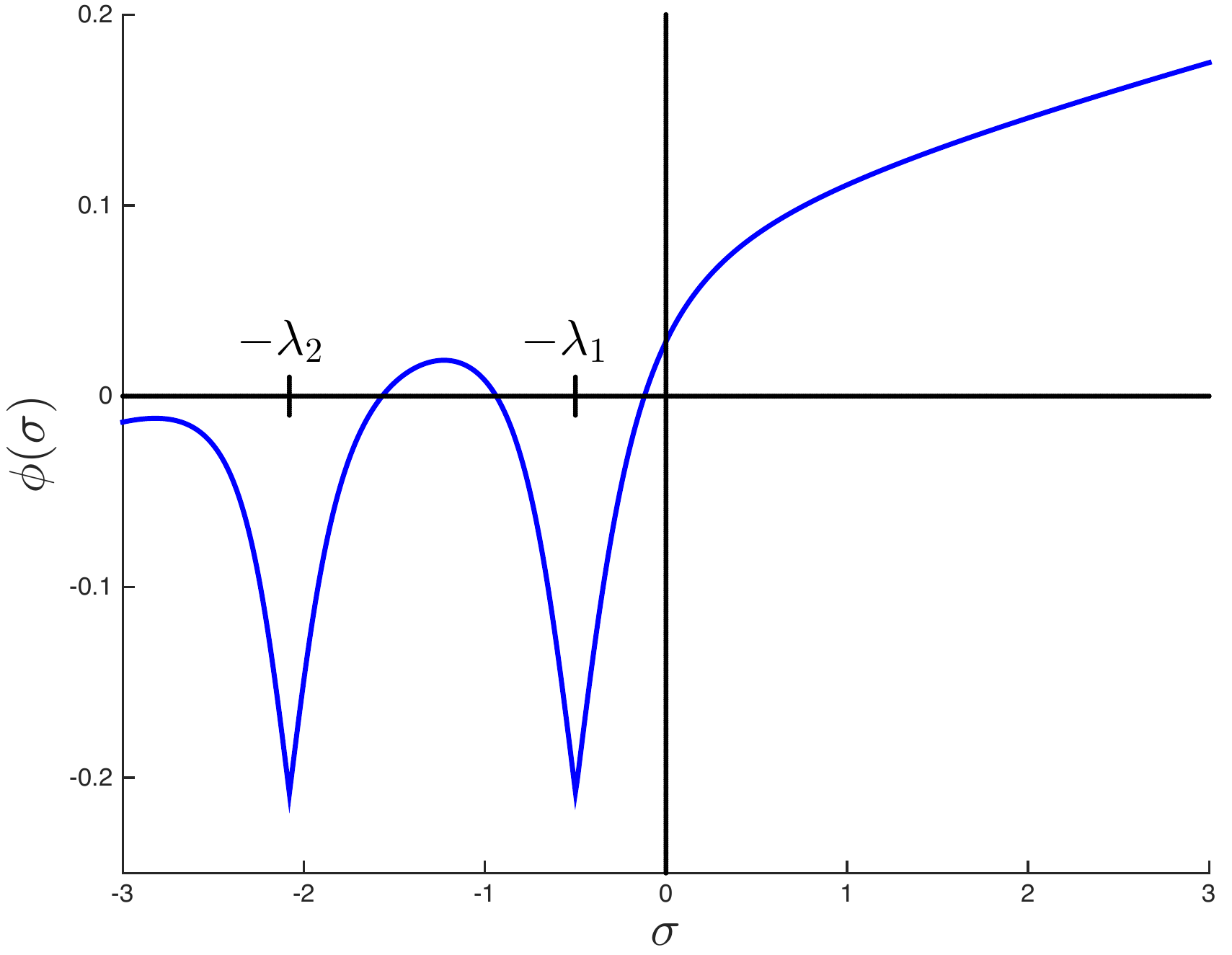}&
	\includegraphics[width=4.2cm]{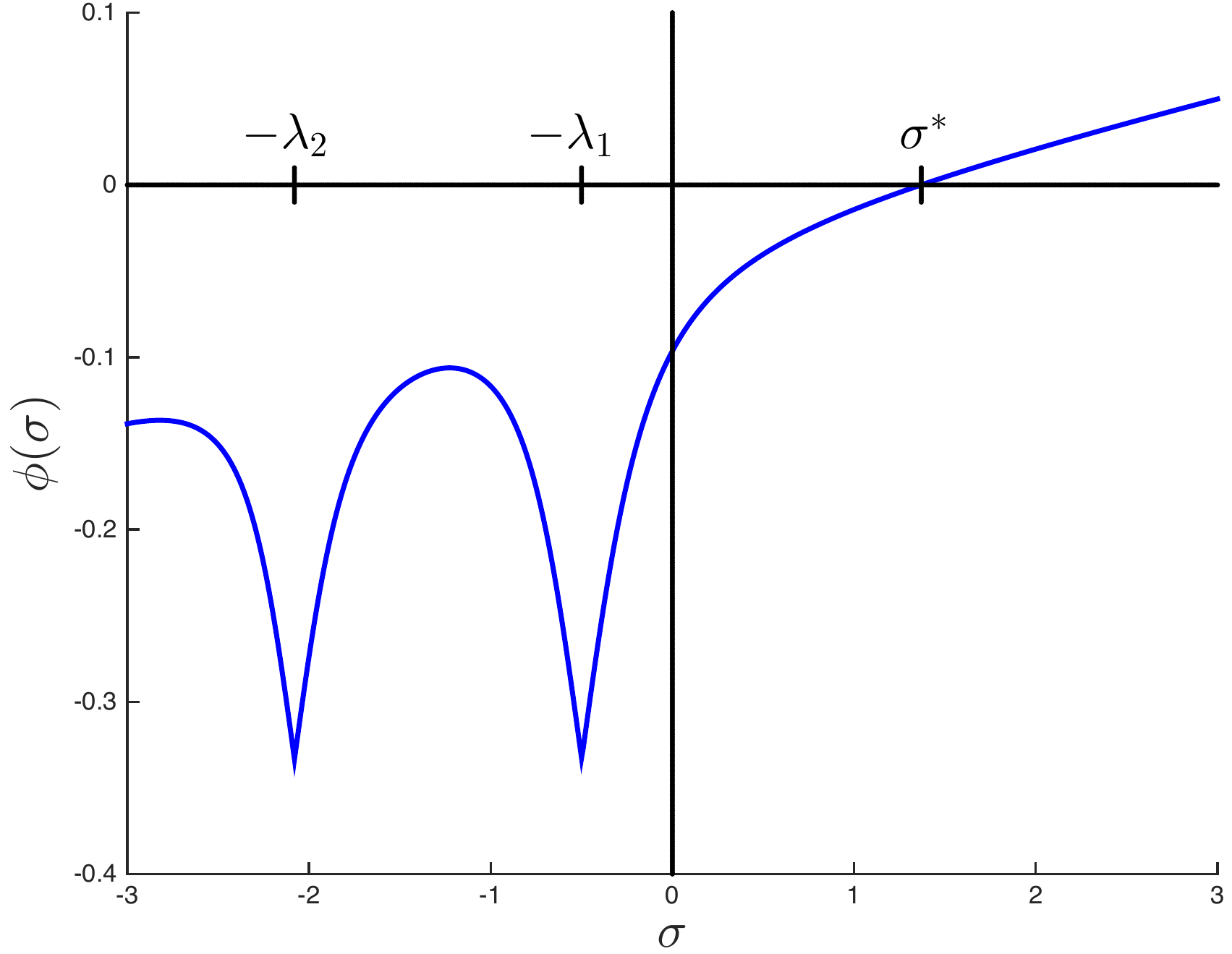}\\
	(a) & (b)
	\end{tabular}
	\caption{Plot of the secular function $\phi(\sigma)$ given in
	\eqref{eqn:secular}. (a) The case when $\phi(0) \ge 0$, which
	implies that the unconstrained minimizer of 
	\eqref{eqn:TrustProblem2} is feasible.
	(b) When $\phi(0) < 0$, there exists $\sigma^* > 0$ such 
	that $\phi(\sigma^*) = 0$, i.e., $v^* = -(\Lambda + \sigma^* I)^{-1}\tilde{g}$
	is well-defined and is feasible.
	}
\end{figure}

The method to obtain $\sigma^*$ is significantly different
that the one used in \cite{BurWX96} in that 
we  explicitly use the eigendecomposition within
Newton's method to compute the optimal $\sigma^*$.  That is, we
differentiate the reciprocal of $\| v(\sigma) \|$ in \eqref{eq:normv}
to compute the derivative of $\phi(\sigma)$ in \eqref{eqn:secular},
obtaining a Newton update that is expressed only in terms of $g_{\parallel}$,
$g_{\perp}$, and the eigenvalues of $B$.  In contrast to 
\cite{BurWX96}, this approach eliminates the
need for matrix solves for each Newton iteration (see Alg.\ 2 in
\cite{BurWX96}).

Given $\sigma^*$ and $v^*$, the optimal $p^*$ is obtained as follows.
Letting $\tau^* = \gamma + \sigma^*$, the solution to the first
optimality condition, $ (B + \sigma^* I )p^* = -g, $ is given by
\begin{eqnarray}
	p^* &=&
	-(B + \sigma^*I) g \nonumber  \\
	&=& -(\gamma I + \Psi M \Psi^T + \sigma^*I)^{-1}g  \nonumber
	\\
	&=& - \frac{1}{\tau^*} 
	\left [
	I - \Psi(\tau^* M^{-1} + \Psi^T\Psi)^{-1}\Psi^T
	\right ]
	g, \label{eqn:p1}
\end{eqnarray}
using the Sherman-Morrison-Woodbury formula.  Algorithm 1 details the
proposed approach for solving the trust-region subproblem.

\begin{algorithm}[ht]
\SetAlgoNoLine
{\bf Algorithm 1:
\LBFGS{} Trust-Region Subproblem Solver}\\
Compute $R$ from the ``thin'' QR factorization of  $\Psi$;\\
Compute the spectral decomposition\\
\quad  \ $RMR^T = V\hat{\Lambda}V^T$ with 
$\hat{\lambda}_1 \le \hat{\lambda}_2 \le \cdots \le \hat{\lambda}_{2(k+1)}$;\\
Let $\Lambda_1 = \hat{\Lambda} + \gamma I$;\\
Define $P_{\parallel} = \Psi R^{-1}V$ and $g_{\parallel} = P_{\parallel}^Tg$;\\
Compute $\| g_{\perp} \|_2 = \sqrt{\| g \|_2^2 - \| g_{\parallel}\|_2^2}$;\\
\textbf{if} $\phi(0) \ge 0$ \textbf{then}\\
\quad $\sigma^* = 0$ and compute $p^*$ from \eqref{eqn:p1} with $\tau^* = \gamma$;\\
\textbf{else} \\
\quad Use Newton's method to find $\sigma^*$;\\
\quad Compute $p^*$ from \eqref{eqn:p1} with $\tau^* = \gamma + \sigma^*$;\\
\textbf{end}
\end{algorithm}

\noindent Algorithm 2 outlines our overall limited-memory \LBFGS{} trust-region approach. 

\begin{algorithm}[ht]
\SetAlgoNoLine
{\bf Algorithm 2: TrustSpa: Limited-Memory \BFGS{}  Trust-Region Method for Sparse Relaxation}\\
Define parameters: $m$, $0 < \tau_1 < 0.5$, $0 < \varepsilon$;\\
Initialize $x_0 \in \Re^{n}$ and compute $g_0 = \nabla \Phi(z_0)$;\\
Let $k = 0$;\\
\textbf{while not converged}\\
\quad \textbf{if} $\| g_k \|_2 \le \varepsilon$ \textbf{then done} \\
\quad Use Algorithm 1 to find $p_k$ that solves \eqref{eqn:TrustProblem};\\
\quad Compute $\rho_k = ( f(z_k+p_k) - f(z_k) )/ q_k(p_k)$;\\
\quad Compute $g_{k+1}$ and update $B_{k+1}$;\\
\quad \textbf{if}  $\rho_k \ge \tau_1$ \textbf{then}\\
\qquad $z_{k+1} = z_k + p_k$;\\
\quad \textbf{else}\\
\qquad $z_{k+1} = z_k;$\\
\quad \textbf{end if} \\
\quad Compute trust-region radius $\delta_{k+1}$;\\
\quad $k \leftarrow k+1$;\\
\textbf{end while}
\end{algorithm}

The method described here guarantees that the trust-region subpoblem 
is solved to high accuracy.
Other quasi-Newton trust-region methods for L-BFGS matrices that solve
to high accuracy include \cite{ErwayM14}, which uses a shifted L-BFGS approach,
and \cite{Burdakov15}, which uses a ``shape-changing'' norm in 
\eqref{eqn:TrustProblem}.

%


\section{NUMERICAL EXPERIMENTS}
\begin{figure}[htbp]\label{setup}
\begin{center}
\begin{tabular}{c} 
(a) Truth $f$ ($\tilde{n} = 4096$, number of nonzeros = 160) \\ 
\hspace*{-.6cm}\includegraphics[scale=0.6]{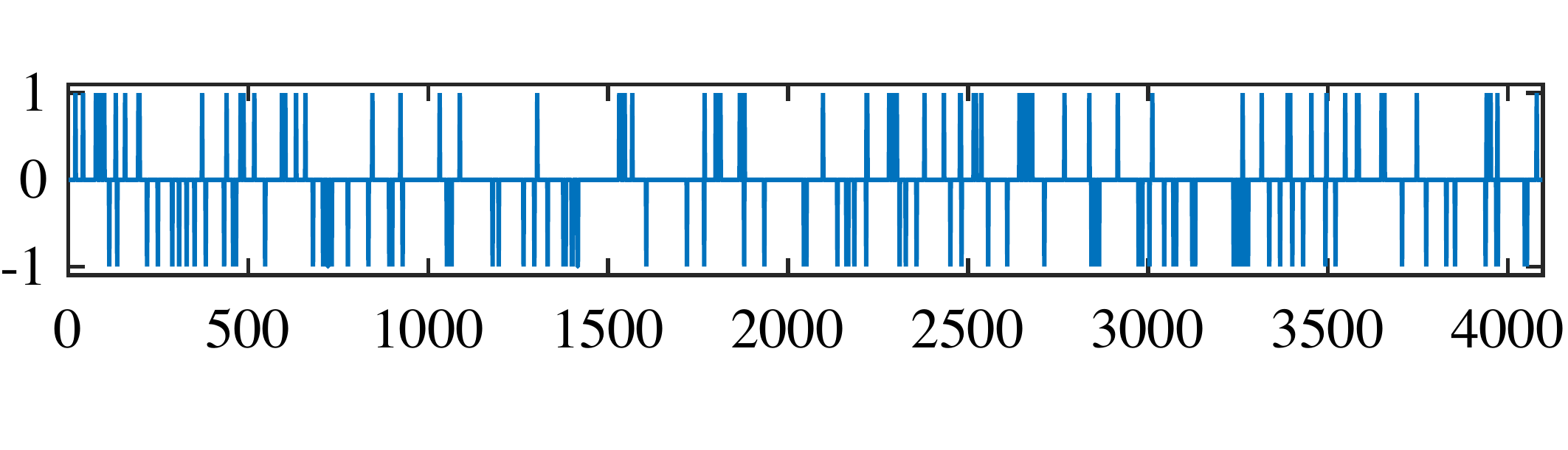} \\
(b) Measurements $y$ ($\tilde{m} = 1024$, noise level = 5\% ) \\ 
\hspace*{-.6cm}\includegraphics[scale=0.6]{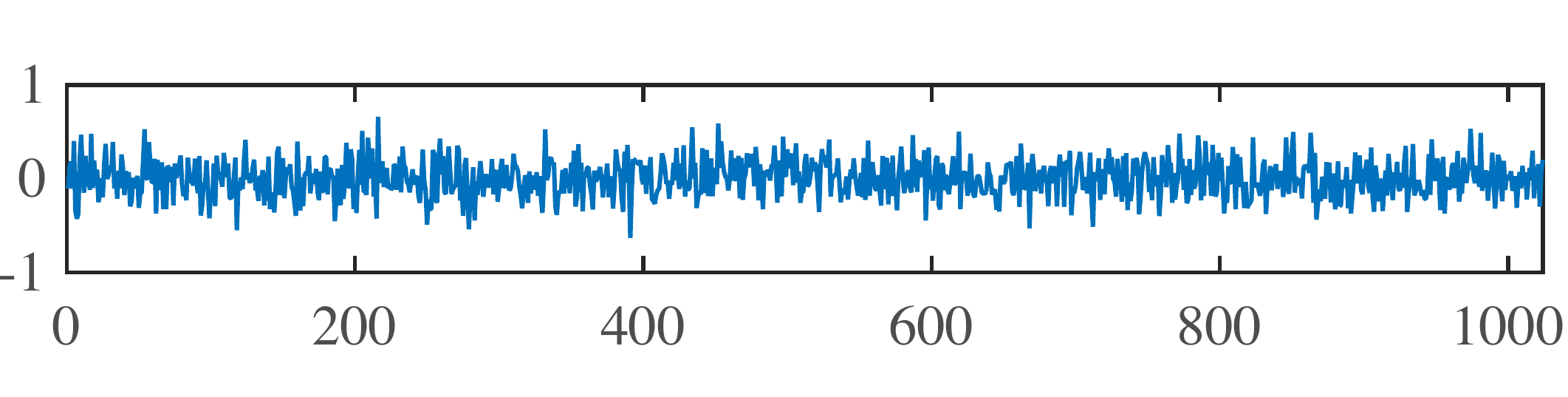} 
\end{tabular}
\caption{Experimental setup: (a) True signal $f$ of size 4,096 with 160 $\pm$ spikes, (b) compressive measurements $y$ 
($\tilde{m} = 1024$) with 5\% Gaussian noise.}
\end{center}
\end{figure}

We call the proposed method, Trust-Region Method for Sparse Relaxation
(TrustSpa Relaxation, or simply TrustSpa).  We evaluate its
effectiveness by reconstructing a sparse signal from Gaussian noise
corrupted low-dimensional measurements. In this experiment, the true
signal $f$ is of size 4,096 with 160 randomly assigned nonzeros with
amplitude $\pm 1$ (see Fig.\ \ref{setup}(a)). We obtain compressive
measurements $y$ of size 1,024 (see Fig.\ \ref{setup}(b)) by
projecting the true signal using a randomly generated system matrix
($A$) from the standard normal distribution with orthonormalized
rows. In particular, the measurements are corrupted by 5\% of Gaussian
noise.

We implemented TrustSpa in MATLAB R2015a using a PC with Intel Core i7 2.8GHz processor with 16GB memory. 
We compared the performance of TrustSpa with the Gradient Projection for Sparse Reconstruction (GPSR) method \cite{Figueiredo2007}
with the Barzilai and Borwein (BB) approach \cite{BARZILAI1988} and without the debiasing option.
Both TrustSpa and GPSR-BB methods are initialized at the same starting point, i.e., zero and terminate 
if the relative objective values do not significantly change, 
i.e, $|\Phi(x^{k+1}) - \Phi(x^{k})|/ |\Phi(x^{k})| \leq 10^{-8}$. The regularization parameter $\tau$ 
in (\ref{eqn:LASSO}) is optimized independently for each algorithm to minimize the 
mean-squared error (MSE = $\tfrac{1}{n}\|\hat{f} - f \|_2^2$, where $\hat{f}$ is an estimate of $f$). 


\begin{figure}[h]
\begin{center}
\begin{tabular}{c} 
 (a) GPSR-BB reconstruction $\hat{f}_{\text{GPSR}}$ (MSE = 1.624e-04) \\ 
\hspace*{-.6cm}\includegraphics[scale=0.6]{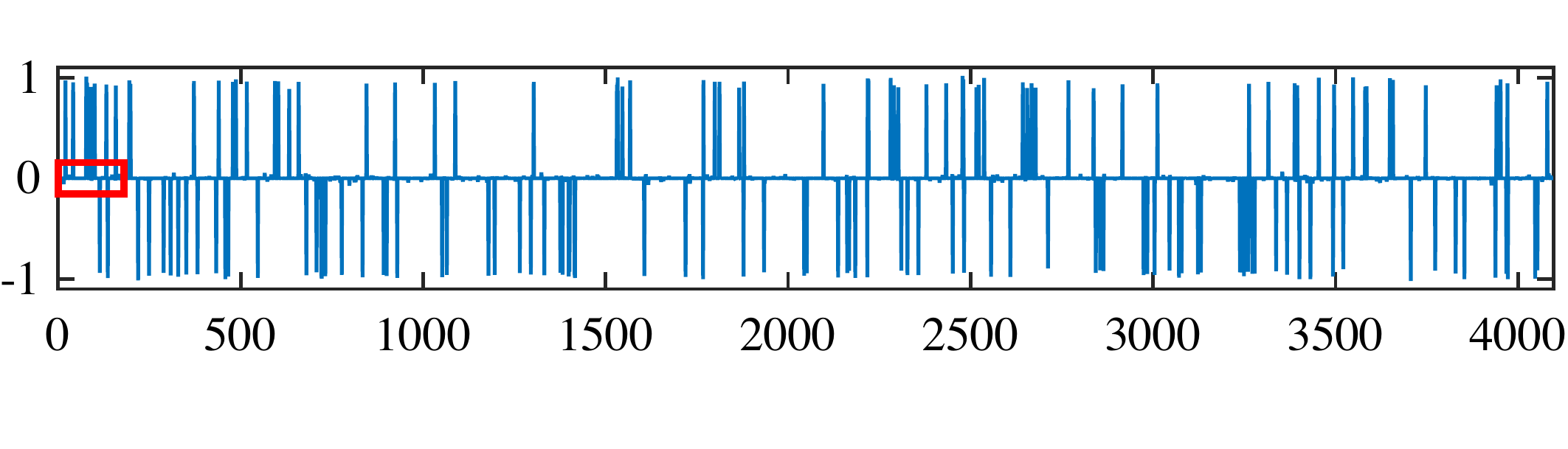} \\
 (b) TrustSpa reconstruction $\hat{f}_{\text{TS}}$ (MSE = 9.347e-05) \\ 
\hspace*{-.6cm} \includegraphics[scale=0.6]{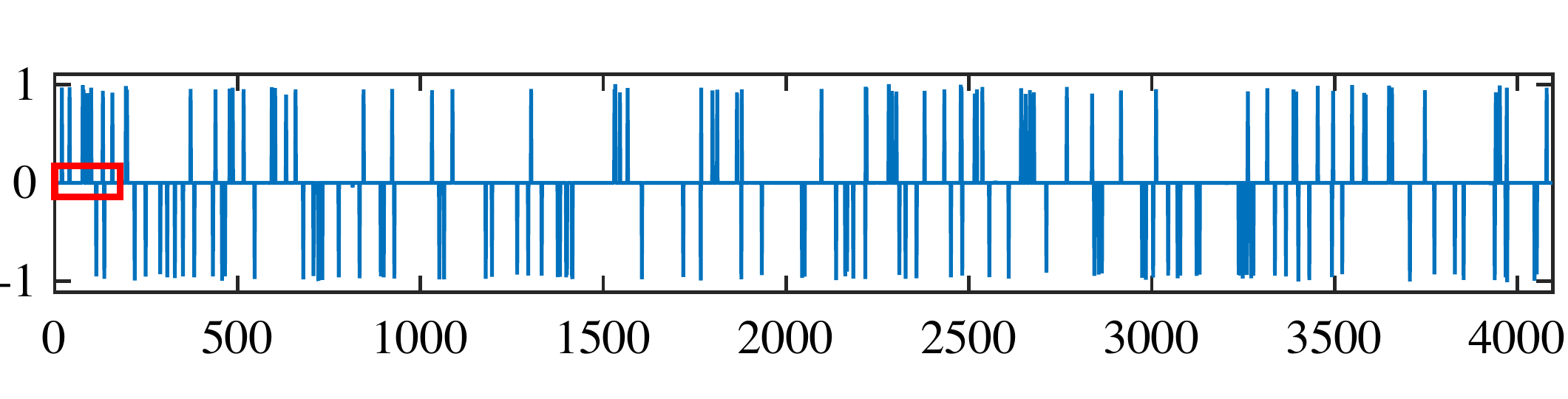} 
\end{tabular}
\caption{(a) GPSR-BB reconstruction, $\hat{f}_{\text{GPSR}}$, (b) TrustSpa reconstruction, $\hat{f}_{\text{TS}}$. 
MSE = $(1/n) \|\hat{f} - f \|_2^2$. Note the lower MSE for the proposed method.}
\label{recons}
\end{center}
\end{figure}

\bigskip

\begin{figure}[h!]
\begin{center}
\begin{tabular}{c} 
(a) Zoomed region of $\hat{f}_{\text{GPSR}}$ \\ 
 \hspace{-0.6cm} \includegraphics[scale=0.6]{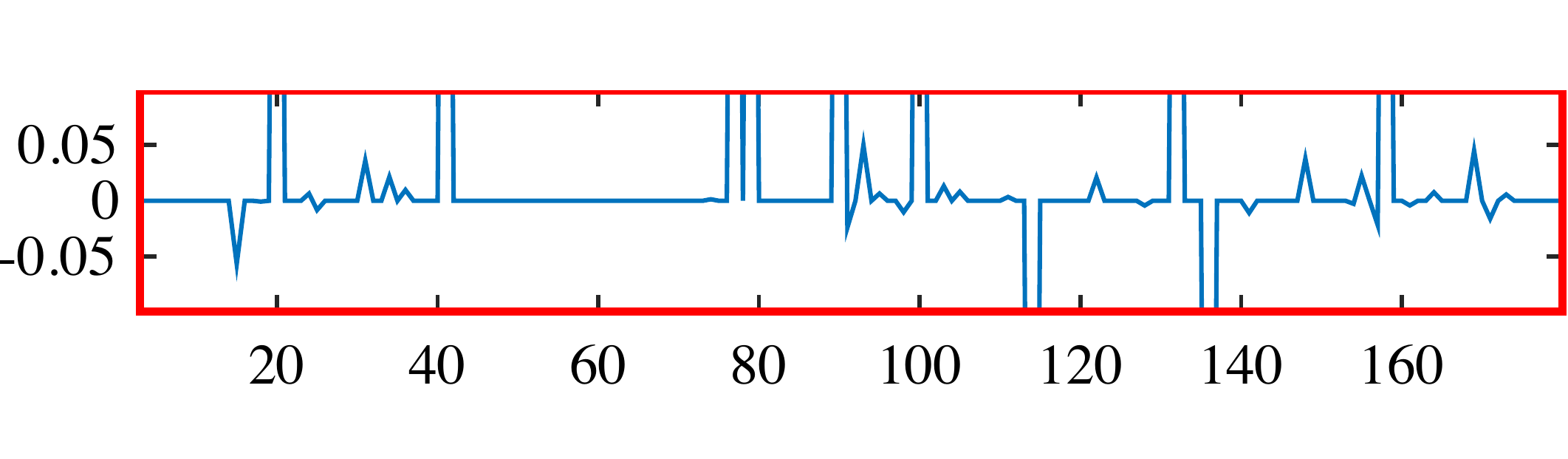} \\
(b)  Zoomed region of  $\hat{f}_{\text{TS}}$  \\ 
 \hspace{-0.6cm} \includegraphics[scale=0.6]{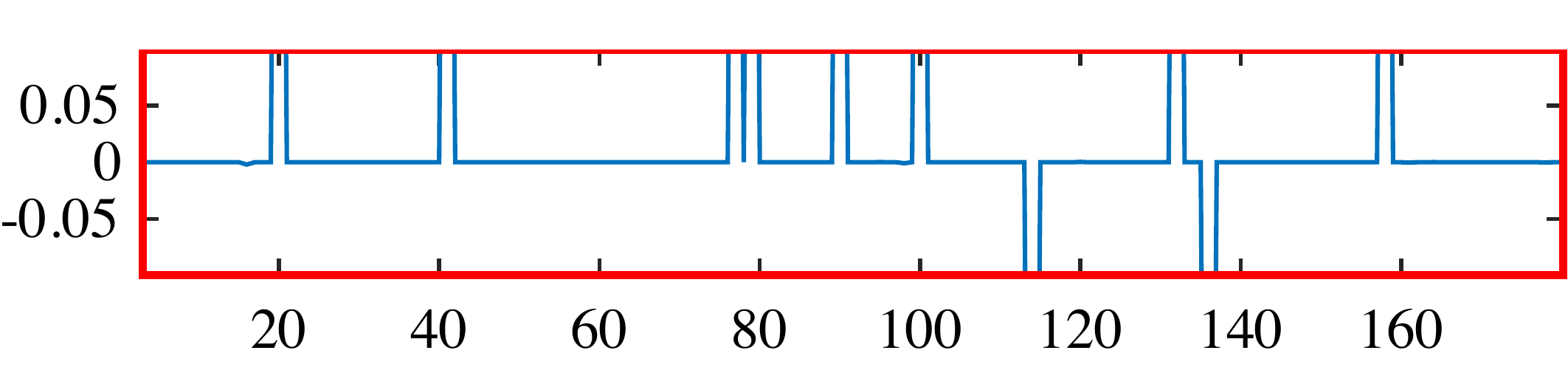} 
\end{tabular}
\caption{Zoomed red-boxed regions in the reconstructions: (a) A zoomed region of $\hat{f}_{\text{GPSR}}$, 
(b) a zoomed region of $\hat{f}_{\text{TS}}$.  Note the presence of artifacts 
in the GPSR-BB reconstruction that are absent in the proposed method's reconstruction.}
\label{zoom}
\end{center}
\end{figure}

\noindent \textbf{Analysis.}  We ran the experiment 10 times with 10
different Gaussian noise realizations.  The average MSE for GPSR-BB for the
10 trials is $1.758 \times 10^{-4}$ and the average computational time is
4.45 seconds. In comparison, the average MSE for TrustSpa is $9.827\times
10^{-5}$, and the average computational time is 3.52 seconds.  For one
particular trial, the GPSR-BB reconstruction, $\hat{f}_{\text{GPSR}}$ (see
Fig.\ \ref{recons}(a)), has MSE $1.624 \times 10^{-4}$ while the TrustSpa
reconstruction, $\hat{f}_{\text{TS}}$ (see Fig.\ \ref{recons}(b)), has MSE
$9.347 \times 10^{-5}$.  Note that the $\hat{f}_{\text{TS}}$ has fewer
reconstruction artifacts (see Fig.\ \ref{zoom}).  Quantitatively,
$\hat{f}_{\text{GPSR}}$ has 786 nonzeros, where the spurious solutions are
between the order of $10^{-2}$ and $10^{-3}$.  In contrast, 
because of the variable transformations used by TrustSpa, the algorithm
terminates with no zero components in its solution; however, only
579 components are greater than $10^{-6}$ in absolute value.
This has the effect of rendering most spurious solutions less visible.

\section{CONCLUSION}

In this paper, we proposed an approach 
for solving the $\ell_2$-$\ell_1$ minimization problem
that arises in compressed sensing and sparse recovery problems.
Unlike gradient projection-type methods like GPSR, which
uses only the current gradient, our approach uses gradients from previous
iterations to obtain a more accurate Hessian approximation.
Numerical experiments show that our proposed approach mitigates
spurious solutions more effectively with a lower average MSE
in a smaller amount of time.

\label{sec:discussion}

\clearpage
\bibliographystyle{IEEEbib}
\bibliography{ErwayMarcia}

\end{document}